\documentclass[12pt]{amsart}
\usepackage{fullpage,url,amssymb,amsmath,graphicx}
\newcommand{\defi}[1]{\textsf{#1}}
\newtheorem{theorem}{Theorem}[section]
\newtheorem{lemma}[theorem]{Lemma}

\theoremstyle{definition}

\newtheorem{example}[theorem]{Example}
\theoremstyle{remark}
\newtheorem{remark}[theorem]{Remark}

\usepackage[backref]{hyperref}
\usepackage[alphabetic,backrefs,lite]{amsrefs}
\newcommand{\D}{\Delta}

\begin{document}

\title{Big rational surfaces}
\author{Damiano Testa, Anthony V\'arilly-Alvarado, Mauricio Velasco}
\address{Mathematical Institute, 24-29 St Giles,
Oxford OX1 3LB, United Kingdom}
\email{adomani@gmail.com}
\address{Department of Mathematics, University of California, 
	Berkeley, CA 94720, USA}
\email{varilly@math.berkeley.edu, velasco@math.berkeley.edu}
\thanks{The first author was partially supported by Jacobs University Bremen, DFG grant STO-299/4-1 and EPSRC grant number EP/F060661/1; the third author is partially supported by NSF grant DMS-0802851.}

\begin{abstract}
We prove that the Cox ring of a smooth rational surface with big anticanonical class is finitely generated.
We classify surfaces of this type that are blow-ups of $\mathbb{P}^2$ at distinct points lying on a (possibly reducible) cubic.
\end{abstract}

\maketitle

\section{Introduction}

\subsection{Mori dream spaces}

Let $X$ be a smooth projective variety over an algebraically closed field. Assume that the Picard group ${\rm Pic}(X)$ is freely generated by the classes of divisors $D_1, D_2,\ldots, D_r$.  The \defi{Cox ring}, or \defi{total coordinate ring}, of $X$ with respect to this choice is given by
\[
{\rm Cox}(X) := \bigoplus_{(m_1,\ldots,m_r)\in\mathbb{Z}^r} {\rm H}^{0}\big( X, \mathcal{O}_X(m_1 D_1+\dots+m_r D_r) \big) ,
\]
with multiplication induced by product of functions in the function field of $X$.

The Cox rings of certain classes of varieties are particularly simple. In the case of toric varieties, for instance, the Cox ring is the ring of polynomial functions on an affine space $\mathbb{A}^d$ with coordinates indexed by the torus-invariant divisors \cite{Cox}*{Theorem~2.1}. Moreover, the variety $X$ can be recovered as a quotient of an open subset of $\mathbb{A}^d = {\rm Spec}({\rm Cox}(X))$ by the action of a torus. More generally, any smooth projective variety $X$ with a \emph{finitely generated} Cox ring can be described in this way:  there is an open subset $\mathcal{T}$ of  ${\rm Spec}({\rm Cox}(X))$ with a canonical torus action and the quotient of $\mathcal{T}$ by this action is isomorphic to $X$. The space $\mathcal{T}$ is an example of a \defi{universal torsor} (see~\cite{CTSII} for a foundational treatment of universal torsors).

Varieties $X$ with finitely generated Cox ring are also distinguished amongst all varieties: the minimal model program on $X$ can be carried out for any divisor.  This privileged position has earned such varieties the name of \defi{Mori dream spaces} \cite{HK}*{Definition~1.10}.   Determining which varieties are Mori dream spaces remains a difficult problem, even in the case of surfaces.  In this paper we are primarily interested in the case in which $X$ is a smooth rational surface with big anticanonical divisor. We prove the following result.

\setcounter{section}{2} \setcounter{theorem}{8}
\begin{theorem}
Let $X$ be a smooth rational surface such that $-K_X$ is big. Then the Cox ring of $X$ is finitely generated.
\end{theorem}
\setcounter{section}{1} \setcounter{theorem}{0}

This result extends a theorem of Hassett \cite{Ha}*{Theorem~5.7}, which states that a smooth rational surface with big and nef anticanonical divisor is a Mori dream space.  As we were completing this paper, we learned that Chen and Schnell obtained an independent proof of Theorem~\ref{eccolo}~\cite{CS}.

\begin{figure}[h]
\begin{center}
\scalebox {0.60}{
\includegraphics{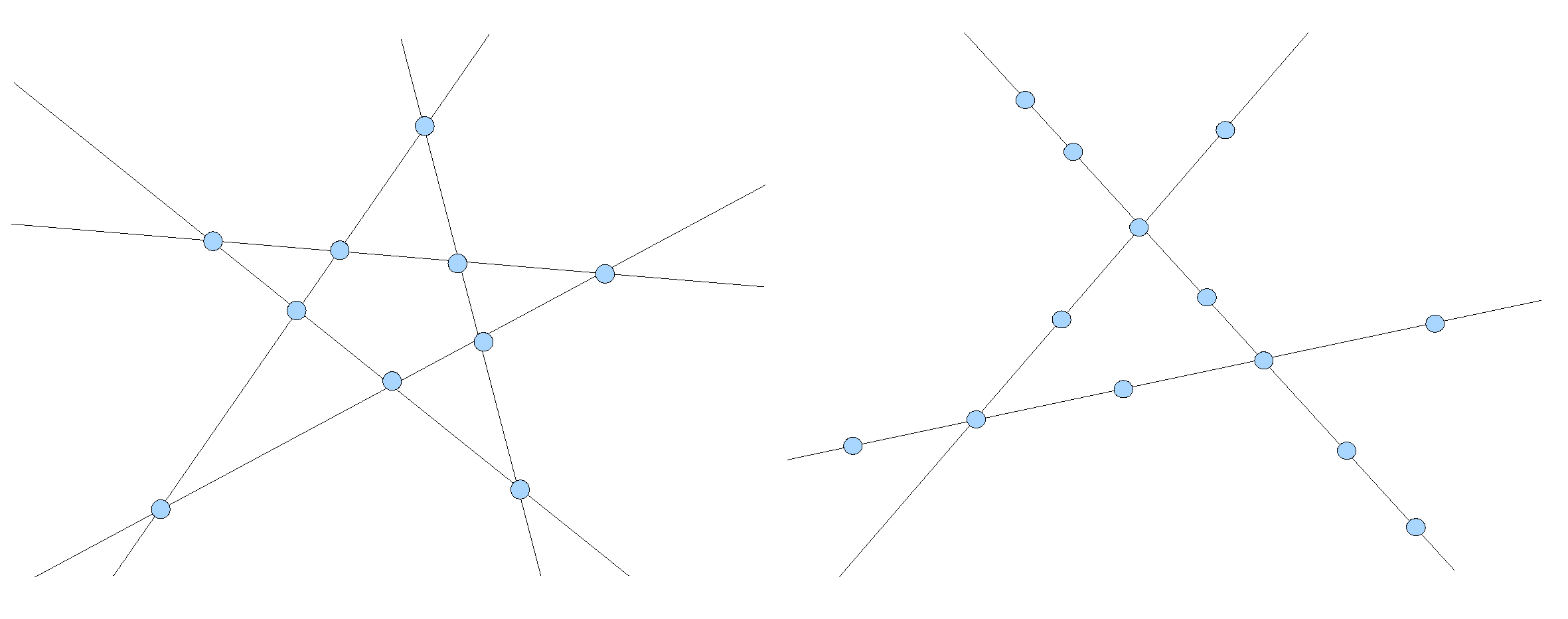}
}
\caption{Examples~\ref{first examples}~(\ref{Castravet}) and~(\ref{E8}). \label{pentagramma}}
\label{F:figura}
\end{center}
\end{figure}

\begin{example}\label{first examples}
We deduce from Theorem~\ref{eccolo} that the surfaces in the following list are Mori dream spaces (see Remark~\ref{first examples explained}).
\renewcommand{\theenumi}{\alph{enumi}}
\begin{enumerate}
\item{Rational surfaces with $K_X^2>0$, or equivalently, rational surfaces with ${\rm rk}({\rm Pic}(X))\leq 9$, are Mori dream spaces.\label{Kx pos}}
\item{Blow-ups of the Hirzebruch surface $\mathbb{F}_n$, $n \geq 1$ at any number of points lying in the union of the negative curve and $n+1$ distinct fibers of the projection.}\label{Hirz's at lots}
\item{Blow-ups of $\mathbb{P}^2$ at $n+1$ points, $n$ of which lie on a (possibly reducible) conic.}\label{P2 on a conic}
\item The surface obtained by blowing up $\mathbb{P}^2$ at the ten points of pairwise intersections of five general lines (Figure~\ref{pentagramma}). Generators for the Cox ring of this surface are determined in~\cite{Castravet}. Harbourne and Ro\'e have also shown that the surface in question is a Mori dream space~(\cite{HQ}).\label{Castravet}
\item{The surface obtained by considering three distinct lines $L_1,L_2,L_3$ in $\mathbb{P}^2$ and blowing up the three pairwise intersections and $2$, $3$ and $5$ additional points on $L_1, L_2$ and $L_3$ respectively.\label{E8}}
\end{enumerate}
\end{example}

To further illustrate the applicability of Theorem~\ref{eccolo}, we classify blow-ups of $\mathbb{P}^2$ at finite sets of points for which $-K_X$ is both big and effective (Theorem~\ref{classi}).  The classification is achieved by associating a root system to each big rational surface, extending a well-known construction for del Pezzo surfaces (Section~\ref{radici}).

There are smooth projective rational surfaces with finitely generated Cox ring, whose anticanonical divisor is not big. For example, by~\cite{Totaro}*{Theorem~5.2} the surface $X$ obtained by blowing up the nine inflection points of a smooth plane cubic has finitely generated Cox ring. However, the anticanonical divisor $-K_X$ is not big since $|{-K_X}|$ contains an irreducible curve and $K_X^2=0$.

\subsection{Previous work} 

Smooth rational surfaces with big anticanonical class have been studied before.  Sakai showed that their anticanonical models have only isolated rational singularities and he provided numerous examples of them (\cite{Sak}).

Harbourne studied the effective and nef cones for rational surfaces with $K_X^2 > 0$ in~\cite{Harbourne}, while~\cite{HGM} proves that the subsemigroup of effective curves is finitely generated for the same class of surfaces. 

The surfaces appearing in Theorem~\ref{classi}~(ii),~(iii) below have been studied by Failla, Lahyane and Molica Bisci in~\cites{FLM,FLM2}; in particular, they investigated the finite generation of the monoid of effective divisor classes modulo algebraic equivalence on these surfaces.

\subsection{Relation to log del Pezzo surfaces}

A log del Pezzo surface (i.e., a Kawamata log terminal pair $(X,\Delta)$ such that $X$ is a normal surface and $-(K_X+\Delta)$ is $\mathbb{Q}$-Cartier and ample) is a Mori dream space \cite{BCHM}*{Corollary~1.3.1}. Clearly, log del Pezzo surfaces have big anticanonical class. It is natural to wonder if the converse is true. It is not. In \S\ref{esempi} we give a family of examples, suggested to us by Chenyang Xu, which shows that the class of smooth rational surfaces with big anticanonical class is strictly larger than the class of log del Pezzo surfaces.

\subsection*{Acknowledgements}
We thank Ana-Maria Castravet, Johan de Jong, David Eisenbud, Brian Harbourne, Brendan Hassett, Se\'an Keel, Bjorn Poonen, Burt Totaro and Chenyang Xu for helpful conversations during the completion of this work.

\section{Big rational surfaces are Mori dream spaces}

From now on, unless otherwise specified, $X$ denotes a smooth projective rational surface over an algebraically closed field with big anticanonical divisor $-K_X$.  We note that ${\rm Pic}(X)$ is a free abelian group and that ${\rm rk}({\rm Pic}(X)) + K_X^2 = 10$. Let $N_1(X)$ be the $\mathbb{R}$-vector space of numerical equivalence classes of curves on $X$; let $NE(X)$ be the cone in $N_1(X)$ of non-negative real combinations of classes of curves on $X$ and let $\overline{NE}(X)$ be its closure.  We denote by $NE(X)_{\mathbb{Z}}$ the cone of non-negative integral linear combinations of classes of curves on $X$.  

\begin{theorem} \label{andra}
Let $X$ be a smooth projective rational surface such that $-K_X$ is big.  The cone of effective divisors on $X$ is finitely generated.  In particular, there are only finitely many integral curves with negative square contained on $X$.
\end{theorem}

\begin{proof}
If the dimension of $N_1(X)$ is at most two, then $X$ is either isomorphic to $\mathbb{P}^2$ or to a Hirzebruch surface and the result is clear.  If the dimension of $N_1(X)$ is at least three, then the result follows from~\cite{Nak}*{Proposition~3.3}. The last statement also follows from \cite{Sak}*{Proposition~4.4}.
\end{proof}

\begin{lemma} \label{hodge}
Let $X$ be a smooth projective surface, let $N$ be a big and nef divisor, and let $C \subset X$ be an effective divisor such that $C \cdot N=0$.  If $C_1 , \ldots , C_r$ are distinct 
irreducible components of $C$, then the matrix $(C_i \cdot C_j)_{i,j}$ is negative definite, and $r \leq \dim (N_1(X))-1$.
\end{lemma}

\begin{proof}
Since $N$ is nef we have  $N\cdot C_i=0$ for all $i$. Since $N$ is big and nef we have $N^2>0$ and the Hodge Index Theorem implies that the matrix $(C_i \cdot C_j)_{1 \leq i,j \leq r}$ is negative definite and has therefore rank $r$.  Thus the vectors $[C_1] , \ldots , [C_r] \in N_1(X)$ are independent and the result follows since they are all contained in the hyperplane orthogonal to $N$.
\end{proof}

\begin{lemma} \label{sezioni}
Let $N$ be a nef divisor on $X$ not linearly equivalent to zero.
\begin{enumerate}
\item \label{sislin}
We have $-K_X \cdot N > 0$ and the linear system $|N|$ has dimension at least one.
\item \label{genar}
If $C \subset X$ is an effective divisor such that $C \cdot N=0$, then the arithmetic genus of $C$ is non-positive; in particular every reduced connected component of $C$ has arithmetic genus zero and every integral component of $C$ is a smooth rational curve.
\end{enumerate}
\end{lemma}

\begin{proof}
(\ref{sislin})
Write $-K_X = A + E$, where $A$ is an ample $\mathbb{Q}$-divisor and $E$ is an effective $\mathbb{Q}$-divisor.  Since $N$ is nef, it is a limit of ample divisors and in particular it is in the closure of the effective cone.  Because $N$ is not linearly equivalent to zero, Kleiman's ampleness criterion \cite{Kl}*{Proposition~IV.2.2} implies that $A \cdot N > 0$, and hence 
$$ -K_X \cdot N = A \cdot N + E \cdot N \geq A \cdot N > 0 $$
since $N$ is nef and $E$ is effective.  We conclude by applying the Riemann-Roch formula to the divisor $N$, together with Serre duality and the fact that $K_X-N$ is the opposite of a big divisor and is therefore not effective.

\noindent
(\ref{genar})
Since $X$ is rational we have ${\rm H}^1(X,\mathcal{O}_X) = (0)$.  From the exact sequence
$$ 0 \longrightarrow \mathcal{O}_X(-C) \longrightarrow \mathcal{O}_X 
\longrightarrow \mathcal{O}_C \longrightarrow 0 $$
we deduce that ${\rm H}^1(C,\mathcal{O}_C)$ is contained in ${\rm H}^2(X,\mathcal{O}_X(-C))$ and, by Serre duality, we have $\dim {\rm H}^2(X,\mathcal{O}_X(-C)) = \dim {\rm H}^0(X,\mathcal{O}_X(K_X+C))$.  By~(\ref{sislin}) we have $-K_X \cdot N>0$ and by assumption $C \cdot N=0$; thus $(K_X +C) \cdot N<0$ which implies that $K_X +C$ is not effective, since $N$ is nef.  It follows that the arithmetic genus of $C$ is non-positive.
\end{proof}

\begin{remark}
Under the additional assumption that $-K_X$ is effective, it is possible to show that $h^0\big(X,\mathcal{O}_X(N)\big) = (N^2 - K_X\cdot N)/2 + 1$ (see \cite{Harbourne2}*{Theorem~III.1 and Lemma~II.2})
\end{remark}

\begin{lemma} \label{nefbpf}
If $N$ is a nef non big divisor on $X$, then $|N|$ is base point free.
\end{lemma}

\begin{proof}
Since $N$ is nef and not big it follows that $N^2 = 0$.  The result is clear if $N=0$; thus from now on we assume that $N \neq 0$.  By Lemma~\ref{sezioni}~(\ref{sislin}), $N$ is linearly equivalent to an effective divisor, and we therefore reduce to the case in which $N$ is effective.  If $N=N_1+N_2$, where $N_1$, $N_2$ are nef divisors, then both $N_1$ and $N_2$ are not big, since otherwise $N$ would be big, and it suffices to show the result for $N_1$ and $N_2$ separately.  Thus we reduce to the case in which $N$ is not the sum of two non-zero nef divisors, since the nef cone contains no lines.

Write $N = a_1 E_1 + \ldots + a_r E_r$, where $E_1 , \ldots , E_r$ are distinct prime divisors and $a_1 , \ldots , a_r \geq 1$ are integers.  The matrix $M := (E_i \cdot E_j)_{1 \leq i,j \leq r}$ is symmetric and has non-negative off-diagonal entries.  Moreover, if $I \subset \{ 1 , \ldots , r\}$ is a subset such that for all $i \in I$ and all $j \in \{ 1 , \ldots , r\} \setminus I$ we have $E_i \cdot E_j = 0$, then $N_1:=\sum _{i\in I} a_i E_i$ and $N_2:=\sum _{j \notin I} a_j E_j$ are nef divisors whose sum equals $N$.  By our reductions, it follows that one of $N_1$, $N_2$ equals zero and therefore that either $I = \emptyset$, or $I = \{ 1 , \ldots , r \}$. Thus $M$ is an irreducible matrix. 

By the Perron-Frobenius Theorem, the largest eigenvalue $\lambda $ of $M$ has an eigenvector $(b_1,\dots,b_r)^t$ with positive entries.  If $\lambda $ were positive, then the effective divisor $b_1E_1 + \cdots + b_rE_r$ would be nef and big, and thus $N$ would itself be big, contradicting the assumptions. Since $M(a_1,\ldots, a_r)^{t}=0$, the largest eigenvalue of $M$ is zero, and by the Perron-Frobenius Theorem it has multiplicity one. Therefore every effective divisor $E$ supported on $\{E_1,\dots,E_r\}$ satisfies $E^2 \leq 0$, with equality if and only if $E$ is proportional to $N$.  

Write $N=P+F$, where $F$ is the fixed divisor of $|N|$ and $P$ is a divisor such that $|P|$ has no base component; in particular we have $P^2 \geq 0$.  Note that $P \neq 0$; otherwise $F$ would have at least two independent global sections, by Lemma~\ref{sezioni}~(\ref{sislin}).  By the above, it follows that $P$ is proportional to $N$ and therefore $F=0$ and $N=P$. We deduce that $|N|$ has no fixed component. Finally, $|N|$ has no base points since the number of its base points is at most $N^2=0$.
\end{proof}

\begin{lemma} \label{bigbpf}
If $N$ is a big and nef divisor on $X$, then $N$ is semiample.
\end{lemma}

\begin{proof}
Let $C \subset X$ denote the union of the integral curves orthogonal to $N$.  By Lemmas~\ref{hodge} and~\ref{sezioni}\eqref{genar} the divisor $C$ satisfies the hypotheses of Artin's contractability criterion~\cite{art}*{Theorem~2.3} and therefore there exists a normal projective surface $X'$ and a birational morphism $X \to X'$ contracting only the connected components of $C$.  By~\cite{art}*{Corollary~2.6} it follows that $N$ is linearly equivalent to a divisor whose support is disjoint from $C$, and therefore $N$ is the pull-back of a Cartier divisor $N'$ on $X'$.  By the Nakai-Moishezon criterion, the divisor $N'$ is ample and hence it is semiample.  Thus its pull-back $N$ is semiample, as we wanted to show.
\end{proof}

We prove a weaker version of~\cite{GalindoMonserrat}*{Theorem~1}.

\begin{lemma} \label{genefa}
The Cox ring of $X$ is generated by global sections supported on the curves with negative self-intersection and generators of the subring $\bigoplus _{N{\rm{ nef}}} {\rm H}^0 (X , \mathcal{O}_X(N))$.
\end{lemma}

\begin{proof}
Let $\mathcal{G} \subset {\rm Cox} (X)$ be a set containing a non-zero section $s_C$ of ${\rm H}^0 (X, \mathcal{O}_X(C))$ for each integral curve $C$ with negative square and a generating set for $\bigoplus _{N{\text{ nef}}} {\rm H}^0 (X , \mathcal{O}_X(N))$.  Fix an ample divisor $A$ on $X$.  We prove by induction on $n$ that for all divisors $D$ on $X$ with $A \cdot D = n$ the vector space ${\rm H}^0 (X, \mathcal{O}_X(D))$ is generated by monomials in $\mathcal{G}$.  The result is clear if $n \leq 0$, since the only effective divisor $D$ with $A \cdot D \leq 0$ is the divisor $D=0$, and the vector space ${\rm H}^0 (X, \mathcal{O}_X)$ is spanned by the empty product of the monomials in $\mathcal{G}$.  Suppose that $n>0$ and that the result is true for all divisors $D'$ such that $A \cdot D' < n$.  Let $D$ be a divisor on $X$ such that $A \cdot D = n$.  If $D$ is either nef or not effective, then there is nothing to prove; so we reduce to the case in which $D$ is effective and not nef.  Therefore there is an integral curve $C$ such that $D \cdot C < 0$, and hence $C^2<0$ and $C$ is contained in the base locus of $|D|$.  Thus the section $s_C$ divides all the vectors in ${\rm H}^0 (X, \mathcal{O}_X(D))$ and the result follows by the inductive hypothesis applied to the divisor $D-C$.
\end{proof}

\begin{lemma}[{\cite{HK}*{Lemma~2.8}}]
\label{zariski}
Let $X$ be a projective variety and let $A_1,\dots,A_r$ be semiample Cartier divisors on $X$. Then the ring
\[
\bigoplus_{(n_1,\dots,n_r) \in \mathbb{Z}^r} {\rm H^0}\big(X,\mathcal{O}_X(n_1A_1 + \cdots n_rA_r)\big)
\]
is finitely generated.\qed
\end{lemma}

\begin{theorem} \label{eccolo}
Let $X$ be a smooth rational surface such that $-K_X$ is big. Then the Cox ring of $X$ is finitely generated.
\end{theorem}

\begin{proof}
By Theorem~\ref{andra} the nef cone of $X$ is finitely generated.  Lemmas~\ref{nefbpf},~\ref{bigbpf}, and~\ref{zariski} together imply that the ring $\bigoplus _{N} {\rm H}^0 (X , \mathcal{O}_X(N))$, as $N$ ranges through all nef divisors, is finitely generated.  By Theorem~\ref{andra} there are only finitely many curves with negative self-intersection on $X$.  Thus the result follows from Lemma~\ref{genefa}.
\end{proof}

\begin{remark}\label{first examples explained}
We briefly explain why the surfaces of Example~\ref{first examples} have big anticanonical class. This is clear for the surfaces of type (\ref{Kx pos}). For a surface $X$ as in~(\ref{Hirz's at lots}), let $\sigma$ and $F$ denote the inverse images of the negative curve and of a fiber, respectively, and let $\tilde{\sigma}$ and $F_1,\dots,F_{n+1}$ denote the strict transforms of $\sigma$ and the special fibers, respectively. We may write
\[
-nK_X = (\sigma + nF) + \Big((n-1)\sigma + n\tilde{\sigma} + n \sum F_i\Big)
\]
which shows that $-nK_X$ is the sum of a big and an effective divisor, whence $-K_X$ is big. For a surface $X$ as in~(\ref{P2 on a conic}), let $p$ be the point of $\mathbb{P}^2$ not on the conic, let $c$ be the strict transform of the conic and let $\ell_1,\dots,\ell_n$ be the strict transforms of the lines through $p$ and each one of the remaining blown-up points. We may write
\[
-nK_X = 2\ell + \Big(\sum \ell_i + (n-1)c\Big);
\]
which shows that $-nK_X$ is the sum of a big and an effective divisor, whence $-K_X$ is big.
Similarly, for the surface (\ref{Castravet}), we note that $-2K_X$ can be written as $\ell + E$ where $\ell$ is the inverse image of the class of a line in $\mathbb{P}^2$ (which is big) and $E$ is effective. For the surface~(\ref{E8}) the divisor $-K_X$ is big by~\S\ref{radici}.
\end{remark}

\section{Big rational surfaces and log del Pezzo surfaces} \label{esempi}

The following family of examples shows that there exist smooth rational surfaces $X$ with big anticanonical divisor which are not log del Pezzo. We show the stronger statement that there is no $\mathbb{Q}$-divisor $\D$ such that $(X,\D)$ is a log canonical pair and $-(K_X+\D)$ is ample.  

Let $n \geq 2$ be an integer and let $h \colon \mathbb{F}_n \to \mathbb{P}^1$ be the Hirzebruch surface with a curve $\bar \sigma $ of square $-n$. Let $k$ be an integer such that $3 \leq k \leq n+1$ and let $a_1 , \ldots , a_k$ be positive integers such that $\sum \frac{1}{a_j} < k-2$.  Choose $k$ distinct integral curves $\bar F_1 , \ldots , \bar F_k \subset \mathbb{F}_n$ contracted by $h$, and for $i \in \{1, \ldots , k\}$ choose $a_i$ distinct points $p^i_1 , \ldots , p^i_{a_i}$ on $\bar F_i \setminus \bar \sigma $.  Let $X$ be the blow-up of $\mathbb{F}_n$ along $\{ p^i_j \mid 1 \leq i \leq k , ~ 1 \leq j \leq a_i\}$; let 
\begin{itemize}
\item $\sigma \subset X$ be the strict transform of the divisor $\bar \sigma $; 
\item $F_i \subset X$ be the strict transform of the divisor $\bar F_i$ for $i \in \{1 , \ldots , k\}$.
\end{itemize}
We have $-K_X = 2 \sigma + (n+2-k) F + \sum F_i$, where $F \subset X$ is the inverse image in $X$ of a fiber of the morphism $h$.  Define $P$ and $N$ as follows
\begin{eqnarray*}
-K_X & = & \overbrace{\frac{n+2-k}{n-\sum \frac{1}{a_j}} \sigma + (n+2-k) F + \sum _i \frac{n+2-k}{a_i (n-\sum \frac{1}{a_j})} F_i} ^{P} + \\[10pt]
& & + \underbrace{\biggl( 2-\frac{n+2-k}{n-\sum \frac{1}{a_j}} \biggr) \sigma + \sum _i \biggl( 1- \frac{n+2-k}{a_i (n-\sum \frac{1}{a_j})} \biggr) F_i} _{N} .
\end{eqnarray*}
Our assumptions on $n, k$ and $\sum \frac{1}{a_j}$ ensure that both $P$ and $N$ are effective. 
Since $P \cdot \sigma = P \cdot F_i = 0$ and $P^2 = \frac{(n+2-k)^2}{n-\sum \frac{1}{a_j}}>0$, it follows that $P$ is big and nef and thus $-K_X$ is big. Additionally $P\cdot N = 0$, so by Lemma~\ref{hodge} the intersection matrix of the support of $N$ is negative definite and therefore $-K_X = P+N$ is the Zariski decomposition of $-K_X$.  

By \cite{Sak}*{Theorem~4.3}, the morphism $f\colon X \to Y$ induced by $|P|$ is a log resolution of $Y := {\rm Proj} \oplus_{m \geq 0}({\rm H}^0\big(X,\mathcal{O}(-mK_X)\big)$ and $f^*(-mK_Y) = mP$ for some $m \gg 0$. We claim that $(Y,0)$ is not a log canonical pair.  Indeed, note that
\[
K_X - f^*(K_Y) = -N 
\]
and that $N$ is supported on divisors contracted by $f$. Thus, the pair $(Y,0)$ is log canonical if and only if $2-\frac{n+2-k}{n-\sum \frac{1}{a_j}} \leq 1$, or equivalently, if and only if $\sum \frac{1}{a_j} \geq k-2$.  

Next, we show that there is no $\mathbb{Q}$-divisor $\D$ such that $(X,\D)$ is a log canonical pair and $-(K_X + \D)$ is ample. To see this, note that if $-(K_X + \D)$ is ample then the divisor
\[
A := -(K_X + \D) - f^*(-K_Y - f_*(\D)),
\]
which is supported on the exceptional locus of $f$, is $f$-ample. By~\cite{Zariski}*{Lemma~7.1}, it follows that all its coefficients are non-positive, and thus $-A$ is effective.  Let $g\colon Z \to X$ be a  log resolution of $(X,\D)$, and let $g^{-1}(\D)$ be the strict transform of $\D$. We have
\smallskip
\[
K_Z + g^{-1}(\D) - g^*(K_X + \D) = \big(K_Z - g^*f^*(K_Y)\big) + \big(g^{-1}(\D) - g^*f^*f_*(\D)\big) - g^*(-A).
\]

\smallskip
\noindent
The coefficients of the divisors $g^{-1}(\D) - g^*f^*f_*(\D)$ and $-g^*(-A)$ are all negative. Since $(Y,0)$ is not log canonical, there is a coefficient of $K_Z - g^*f^*(K_Y)$ which is strictly less than $-1$, and thus $(X,\D)$ is not a log canonical pair.

\section{Blow-ups of the projective plane and root systems} \label{radici}

In this section we classify blow-ups $X$ of $\mathbb{P}^2$ at finite sets of points for which $-K_X$ is both big and effective.  We do so by associating a root system to each big rational surface, extending a well-known construction for del Pezzo surfaces~\cite{M}*{Section~IV.25}.

\begin{lemma} \label{bigma}
Let $X$ be a smooth projective surface and let $\mathcal{D}$ be a set of integral curves on $X$. There is a big divisor whose support is contained in $\mathcal{D}$ if and only if the lattice $\mathcal{D}^{\perp}$ is negative definite.
\end{lemma}

\begin{proof}
\noindent
($\Leftarrow$)  Suppose that $\mathcal{D}^{\perp}$ is negative definite. By the Hodge Index Theorem, there is a divisor $B= \sum _{C \in \mathcal{D}} a_C C$ such that $B^2>0$.  By Riemann-Roch and Serre duality, $h^0(mB)+h^0(K-mB)$ grows (at least) quadratically in $m$ and the same statement holds for $h^0(-mB)+h^0(K+mB)$. Since  $(K-mB)+(K+mB)=2K$ it follows that $h^0(K-mB)$ and $h^0(K+mB)$ cannot both grow quadratically.  We deduce that either $B$ or $-B$ is a big divisor and the result follows.

\noindent
($\Rightarrow$) Suppose that $B= \sum _{C \in \mathcal{D}} a_C C$ is a big divisor with $a_C \in \mathbb{Z}$ for all $C \in \mathcal{D}$. Adding non-negative multiples of the curves in $\mathcal{D}$ to $B$ we reduce to the case in which $B$ is effective; thus the base locus of $B$ is supported on $\mathcal{D}$.  It suffices to show that there exists a big and nef divisor $N$ in the integral span of $\mathcal{D}$, since then $N^2>0$ and $\mathcal{D}^{\perp}\subset N^{\perp}$ is negative definite by the Hodge Index Theorem.  Choose $m \gg 0$ so that the moving part of $mB$ is big. Then, subtracting from $mB$ its base components we obtain the desired divisor $N$.
\end{proof}

\begin{lemma} \label{rootlemma}
Let $\Lambda $ be a negative definite lattice. The set 
$$ R := \bigl\{ \alpha \in \Lambda ~\mid~ \alpha^2 \in \{-1,-2\} \bigr\} $$
is a root system in the span $E$ of $R$.
\end{lemma}

\begin{proof}
We adapt the argument in~\cite{M}*{Section~IV.25}.  To check that $R$ is a root system in $E$ it suffices to verify the axioms in~\cite{Hu}*{Section~III.9}. 

\renewcommand{\theenumi}{R\arabic{enumi}}
\begin{enumerate}
\item {\it The set $R$ is finite, does not contain $0$ and spans $E$}.  This follows from the definition of $E$ and the fact that the pairing is definite.
\item {\it If $\alpha\in R$, then the only multiples of $\alpha$ in $R$ are $\pm\alpha$}. If $m \in \mathbb{R}$ is such that $\alpha , m \alpha \in R$, then $m^2\alpha^2\in \{-1,-2\}\cap \{-m^2,-2m^2\}$ and $m$ is rational since $\alpha, m\alpha \in \Lambda $. We deduce that $m^2=1$.
\item {\it If $\alpha\in R$, then the reflection $\sigma_{\alpha}$ fixing the hyperplane orthogonal to $\alpha$ leaves  $R$ invariant}. The reflection $\sigma_{\alpha}$ is given by 
\[\sigma_{\alpha}(x)=x-2\frac{x\cdot \alpha}{\alpha\cdot\alpha}\alpha.\]
It follows from the definitions that $\sigma_{\alpha} (x)^2=x^2$ for all $x \in \Lambda $.
\item
{\it For every $\alpha, \beta\in R$ we have $\alpha\cdot \beta\in \mathbb{Z}$}.  This property holds for all vectors in $\Lambda $.
\end{enumerate}
The lemma follows.
\end{proof}

Let $X$ be a smooth projective surface and let $\alpha \in N_1(X)_{\mathbb{Z}}$.  It follows from the adjunction formula that $\alpha^2\equiv K_X\cdot \alpha \pmod 2$.  In particular, the quadratic form associated to any sublattice of $N_1(X)_{\mathbb{Z}}$ orthogonal to $K_X$ is even.

\begin{theorem} \label{classi}
Let $\pi\colon X\rightarrow \mathbb{P}^2$ be the blow-up of $\mathbb{P}^2$ at $r$ distinct points. Then $-K_X$ is effective and big if and only if one of the following holds:
\begin{enumerate}
\item $r \leq 8$;\label{del pezzo-ish}
\item a general element of $|{-K_X}|$ consists of the strict transform of a line and a conic where exactly $a$ of the blown-up points lie exclusively on the line, exactly $b$ of the blown-up points lie exclusively on the conic, and either $ab=0$ or $\frac{1}{a}+\frac{4}{b}>1$;\label{break}
\item a general element of $|{-K_X}|$ consists of the strict transform of three lines $L_1$, $L_2$ and $L_3$ where for $i \in \{1,2,3\}$ exactly $a_i$ blown-up points lie exclusively on the line $L_i$, and either $a_1a_2a_3=0$ or $\frac{1}{a_1}+\frac{1}{a_2}+\frac{1}{a_3}>1$.\label{cool case}
\end{enumerate} 
\end{theorem}

\begin{proof}
If $|{-K_X}|$ contains an irreducible divisor $D$, then $D^2 = 9-r$ and, by Lemma~\ref{bigma}, $D$ is big if and only if $r \leq 8$.  Thus we reduce to the case in which every element of $|{-K_X}|$ is reducible and the set of blown-up points  $\mathcal{P}$ is contained in the union of a line and a (possibly reducible) conic.

Suppose that $\mathcal{P}$ is contained in the union of a line $L \subset \mathbb{P}^2$ and an integral conic $C \subset \mathbb{P}^2$.  Let $a$ be the number of points of $\mathcal{P}$ contained in $L \setminus C$ and let $b$ be the number of points of $\mathcal{P}$ contained in $C \setminus L$.  
If $ab=0$, then the result follows from Example~\ref{first examples}~(\ref{P2 on a conic}); thus we reduce to the case $a,b \geq 1$.  Let $\ell \in {\rm Pic}(X)$ be the class of the inverse image of a line, let $e_1 , \ldots , e_a \in {\rm Pic}(X)$ be the classes of the exceptional curves lying above the points of $\mathcal{P}$ in $L\setminus C$ and let $f_1 , \ldots , f_b \in {\rm Pic}(X)$ be the classes of the exceptional curves lying above the points of $\mathcal{P}$ in $C\setminus L$.  The divisor classes $e_1-e_2, e_2-e_3, \ldots , e_{a-1}-e_a$, $f_1-f_2, f_2-f_3, \ldots , f_{b-1}-f_b$, are orthogonal to the components of a general element of $|{-K_X}|$ and are positive roots of a root lattice of type $A_{a-1}(-1) \oplus A_{b-1}(-1)$.  Therefore the intersection form restricted to the span of the above roots is negative definite. On the other hand, the vector $v := ab \ell - b \sum e_i - 2a \sum f_j$ is orthogonal to the components of $-K_X$ and to the root lattice $A_{a-1}(-1) \oplus A_{b-1}(-1)$. It follows from Lemma~\ref{bigma} that $-K_X$ is big if and only if $v^2<0$; hence $-K_X$ is big if and only if $\frac{1}{a}+\frac{4}{b}>1$ and we conclude.

Suppose that $\mathcal{P}$ is contained in the union of three lines $L_1,L_2,L_3 \subset \mathbb{P}^2$.  For $i \in \{1,2,3\}$ let $a_i$ be the number of points of $\mathcal{P}$ contained in $L_i$ and not $L_j$ for $j \neq i$, and let $e^i_1 , \ldots , e^i_{a_i} \in {\rm Pic}(X)$ be the classes of the exceptional curves lying above such points of $\mathcal{P}$; let also $\ell \in {\rm Pic}(X)$ be the class of the inverse image of a line.  If $a_1 a_2 a_3 = 0$, then the result follows from Example~\ref{first examples}~(\ref{P2 on a conic}); thus we reduce to the case $a_1,a_2,a_3 \geq 1$.  The divisor classes $\{e^i_j-e^i_{j+1} ~\mid~ i \in \{1,2,3\} , ~ j \in \{1,\ldots,a_i-1\}\}$ are orthogonal to the components of the element of $|{-K_X}|$ whose image in $\mathbb{P}^2$ is $L_1 + L_2 + L_3$ and are positive roots of a root lattice of type $A_{a_1-1}(-1) \oplus A_{a_2-1}(-1) \oplus A_{a_3-1}(-1)$.  Therefore the intersection form restricted to the span of the above roots is negative definite.  The vector $v := a_1a_2a_3\ell - a_2a_3 \sum e^1_i - a_1a_3 \sum e^2_j - a_1a_2 \sum e^3_k$ is orthogonal to the components of $-K_X$ and to the root lattice $A_{a_1-1}(-1) \oplus A_{a_2-1}(-1) \oplus A_{a_3-1}(-1)$.  It follows from Lemma~\ref{bigma} that $-K_X$ is big if and only if $v^2<0$; hence $-K_X$ is big if and only if $\frac{1}{a_1}+\frac{1}{a_2}+\frac{1}{a_3}>1$, and we conclude.
\end{proof}

\begin{remark}
We describe explicitly the root system of Lemma~\ref{rootlemma} contained in the lattice orthogonal to the components of an element of $|{-K_X}|$, when $X$ is one of the surfaces of Theorem~\ref{classi}.  

For surfaces of type~(\ref{del pezzo-ish}) we recover a subsystem of the usual root system associated to a del Pezzo surface \cite{M}*{Section~IV.25}.

For surfaces of type~(\ref{break}), let $R_{a,b}(-1)$ be the orthogonal complement of the irreducible components of a reducible section of $-K_X$ of the kind mentioned in Theorem~\ref{classi}.  The lattice $R_{a,b}$ is spanned by a root system; a set of positive roots for $R_{a,b}$ is given by 
$$ \begin{array} {ll}
\varepsilon _i := e_i - e_{i+1} & {\text{for }i \in \{1 , 2 , \ldots , a-1\}} , \\[5pt] 
\varphi _j := f_j - f_{j+1} & {\text{for }j \in \{1 , 2 , \ldots , b-1\}} , \\[5pt]
\varepsilon := \ell - e_a - f_1 - f_2 & {\text{if }} a \geq 1 {\text{ and }} b \geq 2.
\end{array} $$
The associated Coxeter graph appears in Figure~\ref{ade1}.

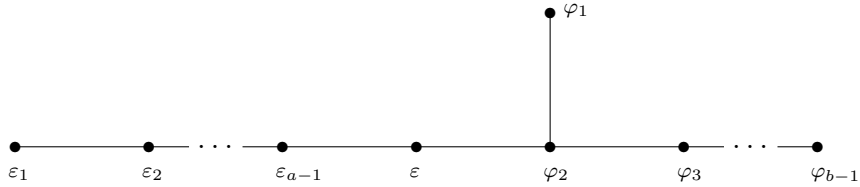
\begin{figure}[h]
\hspace{-115pt} 
\setlength{\unitlength}{0.7in}
\begin{picture}(4.1,1.2)
            \put(0,0.1){\circle*{0.075}}
            \put(1,0.1){\circle*{0.075}}
            \put(2,0.1){\circle*{0.075}}
            \put(3,0.1){\circle*{0.075}}
            \put(4,1.1){\circle*{0.075}}
            \put(4,0.1){\circle*{0.075}}
            \put(5,0.1){\circle*{0.075}}
            \put(6,0.1){\circle*{0.075}}
            \put(0,0.1){\line(1,0){1}}
            \put(1,0.1){\line(1,0){.3}}
            \put(1.7,0.1){\line(1,0){.3}}
            \put(1.36,0.09){\dots}
            \put(2,0.1){\line(1,0){2}}
            \put(4,0.1){\line(0,1){1}}
            \put(5,0.1){\line(1,0){.3}}
            \put(5.7,0.1){\line(1,0){.3}}
            \put(5.36,0.09){\dots}
            \put(4,0.1){\line(1,0){1}}
            \put(-0.05,-0.125){\tiny $\varepsilon _1$}
            \put(0.95,-0.125){\tiny $\varepsilon _2$}
            \put(1.95,-0.125){\tiny $\varepsilon _{a-1}$}
            \put(2.95,-0.125){\tiny $\varepsilon $}
            \put(4.1,1.1){\tiny $\varphi _1$}
            \put(3.95,-0.125){\tiny $\varphi _2$}
            \put(4.95,-0.125){\tiny $\varphi _3$}
            \put(5.95,-0.125){\tiny $\varphi _{b-1}$}
            \end{picture} \caption{The Coxeter graph for the root lattice $R_{a,b}$.} \label{ade1}
            \end{figure}
\noindent
The type of the root lattice $R_{a,b}$ varies with $a$ and $b$: the following are the possibilities.
\begin{itemize}
\item If $ab = 0$, then $R_{a,b}=A_{a+b-1}$.
\item If $b = 2$, then $R_{a,b}=A_a+A_1$.
\item If $b = 3$, then $R_{a,b}=A_{a+2}$.
\item If $b = 4$, or $a=1$ and $b \geq 4$, then $R_{a,b}=D_{a+b-1}$.
\item If $a=2$ and $b = 5$, then $R_{a,b}=E_6$.
\item If $a=3$ and $b = 5$ or $a=2$ and $b = 6$, then $R_{a,b}=E_7$.
\item If $a=4$ and $b = 5$ or $a=2$ and $b = 7$, then $R_{a,b}=E_8$.
\end{itemize}

Similarly, for surfaces of type~(\ref{cool case}), let
$R_{a_1,a_2,a_3}$ be the opposite of the orthogonal complement of the
components of $-K_X$.  The lattice $R_{a_1,a_2,a_3}$ is a root system
of type $A_m+A_n$, $D_n$, $E_6$, $E_7$ or $E_8$, depending on the
values of $a_1$, $a_2$ and $a_3$. A set of positive roots for the root
system is given by $\{\varepsilon ^i_j := e^i_j-e^i_{j+1} ~\mid~ i \in
\{1,2,3\} , ~ j \in \{1,\ldots,a_i-1\}\}$, together with $\varepsilon
:= \ell - e^1_{a_1} - e^2_{a_2} - e^3_{a_3}$ if $a_1,a_2,a_3 \geq 1$.

\noindent
Relabeling the indices if necessary we assume that $a_1 \geq a_2 \geq
a_3$; note also that we necessarily have $a_3 \leq 2$ and if $a_3 =
2$, then $a_2 \leq 3$.  The following are the possibilities.

\begin{itemize}
\item If $a_3 = 0$, then $R_{a_1,a_2,a_3} = A_{a_1-1}+A_{a_2-1}$.
\item If $a_3 = 1$, then $R_{a_1,a_2,a_3} = A_{a_1+a_2-1}$.
\item If $a_2 = a_3 = 2$, then $R_{a_1,a_2,a_3} = D_{a_1+2}$.
\item If $a_2 = 3$ and $a_3 = 2$, then $3 \leq a_1 \leq 5$ and
$R_{a_1,a_2,a_3} = E_{a_1+3}$ (Figure~\ref{F:figura}).
\end{itemize}
The associated Coxeter graph appears in Figure~\ref{ade2}.

\vspace{1.0in}
\begin{figure}[h]
\hspace{-84pt}
\parbox[c]{3.1in}{
\setlength{\unitlength}{0.7in}
\begin{picture}(4.1,0.8)
            \put(0,0.1){\circle*{0.075}}
            \put(1,0.1){\circle*{0.075}}
            \put(2,0.1){\circle*{0.075}}
            \put(3,0.1){\circle*{0.075}}
            \put(3,0.8){\circle*{0.075}}
            \put(3,1.5){\circle*{0.075}}
            \put(3,2.2){\circle*{0.075}}
            \put(4,0.1){\circle*{0.075}}
            \put(5,0.1){\circle*{0.075}}
            \put(6,0.1){\circle*{0.075}}
            \put(0,0.1){\line(1,0){1}}
            \put(1.36,0.09){\dots}
            \put(1,0.1){\line(1,0){.3}}
            \put(1.7,0.1){\line(1,0){.3}}
            \put(2,0.1){\line(1,0){2}}
            \put(3,0.1){\line(0,1){0.7}}
            \put(3,0.8){\line(0,1){.2}}
            \put(3,1.3){\line(0,1){.2}}
            \put(2.97,1.06){\vdots}
            \put(3,1.5){\line(0,1){0.7}}
            \put(4,0.1){\line(1,0){.3}}
            \put(4.7,0.1){\line(1,0){.3}}
            \put(4.36,0.09){\dots}
            \put(5,0.1){\line(1,0){1}}
            \put(-0.05,-0.125){\tiny $\varepsilon^1_1$}
            \put(0.95,-0.125){\tiny $\varepsilon^1_2$}
            \put(1.95,-0.125){\tiny $\varepsilon^1_{a_1 - 1}$}
            \put(2.95,-0.125){\tiny $\varepsilon $}
            \put(3.95,-0.125){\tiny $\varepsilon^3_{a_3 - 1}$}
            \put(3.1,0.77){\tiny $\varepsilon^2_{a_2 - 1}$}
            \put(3.1,1.47){\tiny $\varepsilon^2_2$}
            \put(3.1,2.17){\tiny $\varepsilon^2_1$}
            \put(4.95,-0.125){\tiny $\varepsilon^3_2$}
            \put(5.95,-0.125){\tiny $\varepsilon^3_1$}
           \end{picture}}
           \caption{The Coxeter graph for surfaces of type~(\ref{cool case}).} \label{ade2}
           \end{figure}
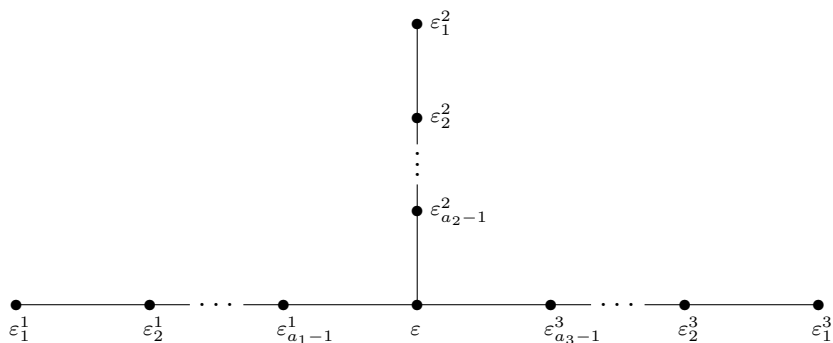
\end{remark}


\begin{bibdiv}
\begin{biblist}

\bib{art}{article}{
   author={Artin, Michael},
   title={Some numerical criteria for contractability of curves on algebraic
   surfaces},
   journal={Amer. J. Math.},
   volume={84},
   date={1962},
   pages={485--496},
   issn={0002-9327},
}

\bib{BCHM}{article}{
    AUTHOR = {Birkar, Caucher},
    author={Cascini, Paolo},
    author={Hacon, Christopher D.},
    author={McKernan, James},    
title={Existence of minimal models for varieties of log general type},
date={2008-8-14},
note={Preprint arXiv:0610203v2 [math.AG]}, 
}

\bib{Castravet}{article}{
    AUTHOR = {Castravet, Ana-Maria},
title={The Cox rings of Log-Fano surfaces (In preparation).},
date={2008},
}

\bib{CS}{article}{
author={Chen, Dawei},
author={Schnell, Christian},
title={Electronic communication},
date={2008-10-24},
}

\bib{CTSII}{article}{
   author={Colliot-Th{\'e}l{\`e}ne, Jean-Louis},
   author={Sansuc, Jean-Jacques},
   title={La descente sur les vari\'et\'es rationnelles. II},
   language={French},
   journal={Duke Math. J.},
   volume={54},
   date={1987},
   number={2},
   pages={375--492},
   issn={0012-7094},
}

\bib{Cox}{article}{
   author={Cox, David A.},
   title={The homogeneous coordinate ring of a toric variety},
   journal={J. Algebraic Geom.},
   volume={4},
   date={1995},
   number={1},
   pages={17--50},
   issn={1056-3911},
}

\bib{FLM}{article}{
   author={Failla, Gioia},
   author={Lahyane, Mustapha},
   author={Molica Bisci, Giovanni},
   title={On the finite generation of the monoid of effective divisor classes on rational surfaces of type $(m,n)$},
   journal={Atti della
Accademia Peloritana dei Pericolanti
Classe di Scienze Fisiche, Matematiche e Naturali},
   volume={LXXXIV},
   date={2006},
   pages={1--9},
   issn={1825 -1242},
}

\bib{FLM2}{article}{
   author={Failla, Gioia},
   author={Lahyane, Mustapha},
   author={Molica Bisci, Giovanni},
   title={The finite generation of the monoid of effective divisor classes
   on Platonic rational surfaces},
   conference={
      title={Singularity theory},
   },
   book={
      publisher={World Sci. Publ., Hackensack, NJ},
   },
   date={2007},
   pages={565--576},
}

\bib{GalindoMonserrat}{article}{
   author={Galindo, Carlos},
   author={Monserrat, Francisco},
   title={The total coordinate ring of a smooth projective surface},
   journal={J. Algebra},
   volume={284},
   date={2005},
   number={1},
   pages={91--101},
   issn={0021-8693},
}

\bib{Harbourne}{article}{
   author={Harbourne, Brian},
   title={Rational surfaces with $K\sp 2>0$},
   journal={Proc. Amer. Math. Soc.},
   volume={124},
   date={1996},
   number={3},
   pages={727--733},
   issn={0002-9939},
}

\bib{Harbourne2}{article}{
   author={Harbourne, Brian},
   title={Anticanonical rational surfaces},
   journal={Trans. Amer. Math. Soc.},
   volume={349},
   date={1997},
   number={3},
   pages={1191--1208},
   issn={0002-9947},
}

\bib{HQ}{article}{
author={Harbourne, Brian},
title={Electronic communication},
date={2009-1-9},
}

\bib{HGM}{article}{
    AUTHOR = {Harbourne, Brian},
    author={Geramita, Anthony V.},
    author={Migliore, Juan},
title={Classifying Hilbert functions of fat point subschemes in $\mathbb{P}^2$},
date={2008-8-8},
note={To appear in \emph{Collect.\  Math.}; math/0803.4113v2}, 
}

\bib{Ha}{article}{
   author={Hassett, Brendan},
   title={Rational surfaces over nonclosed fields},
   journal={to appear in the Proceedings of the 2006 Clay Summer School},
}

\bib{HK}{article}{
   author={Hu, Yi},
   author={Keel, Se{\'a}n},
   title={Mori dream spaces and GIT},
   journal={Michigan Math. J.},
   volume={48},
   date={2000},
   pages={331--348},
   issn={0026-2285},
}

\bib{Hu}{book}{
   author={Humphreys, James E.},
   title={Introduction to Lie algebras and representation theory},
   series={Graduate Texts in Mathematics},
   volume={9},
   publisher={Springer-Verlag},
   place={New York-Berlin},
   date={1978},
   pages={xii+171},
   isbn={0-387-90053-5},
}

\bib{Kl}{article}{
   author={Kleiman, Steven L.},
   title={Toward a numerical theory of ampleness},
   journal={Ann. of Math. (2)},
   volume={84},
   date={1966},
   pages={293--344},
   issn={0003-486X},
}

\bib{M}{book}{
    AUTHOR = {Manin, Yuri I.},
     TITLE = {Cubic forms: algebra, geometry, arithmetic},
 PUBLISHER = {North-Holland Publishing Co.},
   ADDRESS = {Amsterdam},
      YEAR = {1974},
     PAGES = {vii+292},
      ISBN = {0-7204-2456-9},
}

\bib{Nak}{article}{
    AUTHOR = {Nakayama, Noboru},
     TITLE = {Classification of log del Pezzo surfaces of index two},
   journal={J. Math. Sci. Univ. Tokyo},
   volume={14},
   date={2007},
   pages={293--498},
}

\bib{Sak}{article}{
    AUTHOR = {Sakai, Fumio},
     TITLE = {Anticanonical models of rational surfaces},
   journal={Math. Ann.},
   volume={269},
   number={3},
   date={1984},
   pages={389--410},
}

\bib{Totaro}{article}{
   author={Totaro, Burt},
   title={Hilbert's 14th problem over finite fields and a conjecture on the
   cone of curves},
   journal={Compos. Math.},
   volume={144},
   date={2008},
   number={5},
   pages={1176--1198},
   issn={0010-437X},
}


\bib{Zariski}{article}{
   author={Zariski, Oscar},
   title={The theorem of Riemann-Roch for high multiples of an effective
   divisor on an algebraic surface},
   journal={Ann. of Math. (2)},
   volume={76},
   date={1962},
   pages={560--615},
   issn={0003-486X},
}

\end{biblist}
\end{bibdiv}
\end{document}